\documentstyle[12pt]{article}

\def\fr{\frac}
\def\la{\lambda}

\def\Si{\Sigma}

\def\de{\delta}

\def\a{a^{\dagger}}

\newcommand{\bd}{\begin{displaymath}}
\newcommand{\ed}{\end{displaymath}}
\newcommand{\bb}{\begin{equation}}
\newcommand{\ee}{\end{equation}}
\newcommand{\bea}{\begin{eqnarray}}
\newcommand{\eea}{\end{eqnarray}}

\begin{document}
\baselineskip 1.8 \baselineskip


\vspace{.2cm}

\begin{center}
\Large {\bf Two parameter Deformed Multimode Oscillators and  q-Symmetric
States 
}
\\[1cm]
\large W.-S.Chung \\[.3cm]
\normalsize  
Department of Physics and Research Institute of Natural Science,  \\
\normalsize  Gyeongsang National University,   \\
\normalsize   Jinju, 660-701, Korea
\end{center}

\vspace{0.5cm}
\begin{abstract}
Two types of the coherent states for two parameter deformed multimode oscillator
system are investigated. Moreover, two parameter deformed $gl(n)$ algebra and 
deformed symmetric states are constructed.
\end{abstract}

\setcounter{page}{1}
\section{Introduction}

Quantum                                                               
groups                                                                      
or                        
q-deformed                    
Lie         
algebra                
implies   
some    
specific 
deformations 
of                    classical                      Lie           algebras.

From           
a             
mathematical      
point   
of  
view,     
it 
is  
a 
non-commutative                                                  
associative                                                                 
Hopf         
algebra.         
The       
structure        
and 
representation 
theory    of                                   quantum               groups  
have                 
been             
developed        
extensively  
by  
Jimbo    
[1] 
and 
Drinfeld                                                                
[2].                                                                        
            
The                                   
q-deformation                    
of        
Heisenberg       
algebra   
was     
made  
by 
Arik and Coon [3], Macfarlane [4] and Biedenharn [5].
Recently                                                               
there                                               
has                           
been                    
some              
interest               
in    
more     
general   
deformations 
involving                                                                 
an                                                                          
arbitrary                      
real                           
functions     
of           
weight    
generators  
and   
including 
q-deformed algebras as a special case [6-10].

\def\a{a^{\dagger}}

\def\q{q^{-1}}
\def\a{a^{\dagger}}
Recently  Greenberg [11]  has  studied the  following  q-deformation of
multi  mode boson 
algebra:
\bd
a_i \a_j -q \a_j a_i=\de_{ij},
\ed
where the deformation parameter $q$ has to be real.
The main problem of Greenberg's approach is that we can not derive the relation
among $a_i$'s operators at all.
In order to resolve this problem, Mishra and Rajasekaran [12] generalized 
the algebra to complex parameter $q$ with $|q|=1$ and another real deformation
parameter $p$.  
In this paper we use the result  of ref [12] to construct two types of
coherent states and 
q-symmetric ststes.

\def\a{a^{\dagger}}

\def\q{q^{-1}}
\def\a{a^{\dagger}}
\section{Two Parameter Deformed Multimode Oscillators}
\subsection{ Representation and Coherent States}
In this subsection we discuss the algebra given in ref [12] and develop its
reprsentation.
Mishra and Rajasekaran's algebra for multi mode oscillators is given by
\begin{displaymath}
a_i \a_j =q \a_j a_i~~~(i<j)
\end{displaymath}
\begin{displaymath}
a_i\a_i -p\a_i a_i=1
\end{displaymath}
\begin{equation}
a_ia_j=\q a_j a_i ~~~~(i<j),
\end{equation}
where $i,j=1,2,\cdots,n$.
In this case we can say that $\a_i$ is a hermitian adjoint of $a_i$.

The fock space representation of the algebra (1) can be easily constructed
by introducing
the hermitian number operators $\{ N_1, N_2,\cdots, N_n \}$ obeying
\bb
[N_i,a_j]=-\de_{ij}a_j,~~~[N_i,\a_j]=\de_{ij}\a_j,~~~(i,j=1,2,\cdots,n).
\ee
From the second relation of eq.(1) and eq.(2), the relation between the
number operator
and creation and annihilation operator is given by
\bb
\a_ia_i =[N_i]=\fr{p^{N_i}-1}{p-1}
\ee
or
\bb
N_i=\sum_{k=1}^{\infty}\fr{(1-p)^k}{1-p^k}(\a_i)^ka_i^k.
\ee

\def\z{|0,0,\cdots,0>}
Let $\z$ be the unique ground state of this system satisfying
\bb
N_i\z=0,~~~a_i\z=0,~~~(i,j=1,2, \cdots,n)
\ee
\def\n{|n_1,n_2,\cdots,n_n>}
\def\np{|n_1,\cdots,n_i+1.\cdots ,n_n>}
\def\nm{|n_1,\cdots,n_i-1.\cdots ,n_n>}
and $\{\n| n_i=0,1,2,\cdots \}$ be the complete set of the orthonormal 
number eigenstates obeying
\bb
N_i\n=n_i\n
\ee
and
\bb
<n_1,\cdots,   n_n|n^{\prime}_1,\cdots,n^{\prime}_n>=\de_{n_1
n_1^{\prime}}\cdots\de_{n_2 
n_2^{\prime}}.
\ee
If we set
\bb
a_i\n=f_i(n_1,\cdots,n_n) \nm,
\ee
we have, from the fact that $\a_i$ is a hermitian adjoint of $a_i$,
\bb
\a_i\n=f^*(n_1,\cdots, n_i+1, \cdots, n_n) \np.
\ee
Making use  of relation $ a_i  a_{i+1} = \q a_{i+1}  a_i $ we find  the
following relation  
for $f_i$'s:
\bd
q\fr{f_{i+1}(n_1,\cdots, n_n)}{f_{i+1}(n_1, \cdots, n_i-1, \cdots, n_n}
=\fr{f_i(n_1,\cdots, n_n)}{f_i(n_1,\cdots,n_{i+1}-1, \cdots, n_n)}
\ed
\bb
|f_i( n_1, \cdots, n_i+1, \cdots, n_n)|^2 -p |f_i(n_1, \cdots, n_n)|^2=1.
\ee
Solving the above equations we find
\bb
f_i(n_1,\cdots, n_n)=q^{\Si_{k=i+1}^n n_k}\sqrt{[n_i]},
\ee
where $[x]$ is defined as 
\bd
[x]=\fr{p^x-1}{p-1}.
\ed

Thus the representation of this algebra becomes
\bd
a_i|n_1,\cdots, n_n>=q^{\Si_{k=i+1}^n n_k}\sqrt{[n_i]}|n_1,\cdots,
n_i-1,\cdots, n_n>~~~
\ed
\bb
\a_i|n_1,\cdots,   n_n>=q^{-\Si_{k=i+1}^n   n_k}\sqrt{[n_i+1]}|n_1,\cdots,
n_i+1,\cdots, 
n_n>.~~~
\ee

The general eigenstates $\n$ is obtained by applying 
$\a_i$'s operators to the ground state $\z$:
\bb
\n =\fr{(\a_n)^{n_n}\cdots (\a_1)^{n_1} }{\sqrt{[n_n]!\cdots[n_1]!}}\z,
\ee
where 
\bd
[n]!=[n][n-1]\cdots[2][1],~~~[0]!=1.
\ed

The coherent states for $gl_q(n)$ algebra
is usually defined as
\bb
a_i|z_1,\cdots,z_i,\cdots,z_n>_-=z_i|z_1,\cdots,z_{i},\cdots,z_n>_-.
\ee
\def\q{\sqrt{q}}
From the $gl_q(n)$-covariant oscillator algebra we obtain the following
commutation 
relation between $z_i$'s and $z^*_i$'s, where $z^*_i$ is a complex conjugate 
of $z_i$;
\bd
z_iz_j=q z_j z_i,~~~~(i<j),
\ed
\bd
z^*_iz^*_j=\fr{1}{q}z^*_jz_i,~~~~(i<j),
\ed
\bd
z^*_iz_j=q z_j z^*_i,~~~~(i \neq j)
\ed
\bb
z^*_iz_i=z_iz^*_i.
\ee
Using these relations the coherent state becomes

\bb
|z_1,\cdots,z_n>_-=c(z_1,\cdots,z_n)\Si_{n_1,\cdots,n_n=0}^{\infty}
\fr{z_n^{n_n}\cdots z_1^{n_1}}{\sqrt{[n_1]!\cdots[n_n]!}}\n.
\ee
Using the eq.(13) we can rewrite eq.(16) as
\bb
|z_1,\cdots,z_n>_-=c(z_1,\cdots,z_n)e_p(z_n\a_n)\cdots e_p(z_1
\a_1)\z, \ee
where
\bd
e_p(x)=\Si_{n=0}^{\infty}\fr{x^n}{[n]!}
\ed
is a deformed exponential function.

In order to obtain the normalized coherent states, we should impose the
condition
$~{}_<z_1,\cdots,z_n|z_1,\cdots,z_n>_-=1$. Then  the normalized coherent
states are given 
by
\bb
|z_1,\cdots,z_n>_-=\fr{1}{\sqrt{e_p(|z_1|^2)\cdots e_p(|z_n|^2)}}
e_p(z_n\a_n)\cdots e_p(z_1 \a_1)\z,
\ee
where $|z_i|^2=z_iz^*_i=z^*_iz_i$.

\subsection{Positive Energy Coherent States}
The purpose of this subsection is to obtain another type of 
coherent states for algebra (1).
In order to do so , it is convenient to introduce
n subhamiltonians as follows
\begin{displaymath}
H_i=\a_ia_i-\nu,
\end{displaymath}
where
\begin{displaymath}
\nu=\frac{1}{1-p}.
\end{displaymath}
Then the commutation relation between the subhamiltonians and mode operators
are given by
\bb
H_i\a_j=(\de_{ij}(p-1)+1)\a_jH_i,~~~~[H_i,H_j]=0.
\ee
Acting  subhamiltonian on the number eigenstates gives
\bb
H_i\n =-\fr{p^{n_i}}{1-p}\n
\ee
Thus the energy becomes negative when $0<p<1$.
As was noticed in ref [13], for the positive energy states it is not $a_i$
but $\a_i$ that
play a role of the lowering operator:
\def\l{\lambda}
\def\lp{|\l_1p^{n_1},\cdots,\l_n p^{n_n}>}
\def\lpp{|\l_1p^{n_1},\cdots,\l_ip^{n_i+1},\cdots,\l_n p^{n_n}>}
\def\lpm{|\l_1p^{n_1},\cdots,\l_ip^{n_i-1},\cdots,\l_n p^{n_n}>}
\def\lm{|\l_1p^{-n_1},\cdots,\l_n p^{-n_n}>}

\begin{displaymath}
H_i\lp
=\lambda_i p^{n_i}
\lp
\end{displaymath}
\begin{displaymath}
\a_i\lp
=q^{-\Si_{k=i+1}^n n_k}\sqrt{\lambda_i p^{n_i+1}+\nu}\lpp
\end{displaymath}
\bb
a_i\lp
=q^{\Si_{k=i+1}^n n_k}\sqrt{\lambda_i p^{n_i}+\nu}\lpm,
\ee
where $ \lambda_1, \cdots,\l_n >0$.

Due to this fact, it is natural to define coherent states 
corresponding to the representation (21)
 as the eigenstates of $\a_i$'s:

\def\s{\Sigma}
\def\fr{\frac}
\def\l{\lambda}
\def\m{\mu}
\def\z{|z_1,\cdots,z_n>_+}
\def\zz{|z_1,\cdots,z_i,qz_{i+1},\cdots,q z_n>_+}
\bb
\a_i\z =z_i \z
\ee
Because the representation  (21) depends on n  free paprameters $ \la_i$'s
, the coherent 
states
$\z$ can take different forms.

If we assume that the positive energy states are normalizable, 
i.e.$~~~$   $<\la_1   p^{n_1},\cdots   ,\l_n  p^{n_n}|\la_1
p^{n_1^{\prime}},\cdots,\l_n 
p^{n_n^{\prime}}>=\de_{n_1 n_1^{\prime}}\cdots\de_{n_nn_n^{\prime}}$, and
form exactly one 
series 
for some 
fixed $\la_i$'s, then we can obtain
\bd
\z
\ed
\bb
=C \s_{n_1,\cdots,n_n=-\infty}^{\infty}
\left[\Pi_{k=0}^n 
\frac{p^{\fr{n_k(n_k-1)}{4}}}
{\sqrt{(-\fr{\nu}{\l_k};p)_{n_k}}}
\left( \fr{1}{\sqrt{\l_k}}\right)^{n_k}\right]z_n^{n_n}\cdots
z_1^{n_1}\lm. \ee

If we demand that ${}_+<z_1,\cdots,z_n|z_1,\cdots,z_n>_+=1$, we have
\bb
C^{-2} =
\Pi_{k=1}^n{}_0\psi_1(-\fr{\nu}{\l_k};p,-\fr{|z_k|^2}{\l_k})
\ee
where  bilateral p-hypergeometric series ${}_0\psi_1(a;p,x)$is defined by [14]
\bb
{}_0 \psi_1(a
;p ,x) 
=\s_{n=-\infty}^{\infty}
\fr{(-)^n p^{n(n-1)/2}}{(a;p)_{n}}x^n.
\ee

\def\l{\lambda}
\def\m{\mu}

\subsection{Two Parameter Deformed $gl(n)$ Algebra}
The purpose of this subsection is to derive the deformed $gl(n)$ algebra
from the deformed multimode oscillator algebra.
The multimode oscillators given in eq.(1) can be arrayed in bilinears to
construct the generators
\bb
E_{ij}=\a_i a_j.
\ee
From the fact that $\a_i$ is a hermitian adjoint of $a_i$, we know that
\bb
E^{\dagger}_{ij}=E_{ji}.
\ee
Then the deformed $gl(n)$ algebra is obtained from the algebra (1):
\bd
[E_{ii},E_{jj}]=0,
\ed
\bd
[E_{ii},E_{jk}]=0,~~~(i\neq j \neq k )
\ed
\bd
[E_{ij},E_{ji}]=E_{ii}-E_{jj},~~~(i \neq j )
\ed
\bd
E_{ii}E_{ij}-p E_{ij} E_{ii}=E_{ij},~~~(i \neq j)
\ed
\bd
E_{ij}E_{ik}=
\cases{
q^{-1}E_{ik}E_{ij} & if $ j<k$ \cr
qE_{ik}E_{ij} & if $ j>k$ \cr}
\ed

\bb
E_{ij}E_{kl}=q^{2(R(i,k)+R(j,l)-R(j,k)-R(i,l))}E_{kl}E_{ij},~~~(i \neq j
\neq k \neq l),
\ee
where the symbol $ R(i,j)$ is defined by
\bd
R(i,j)=\cases{
1& if $i>j$\cr
0& if $i \leq j $ \cr }
\ed

This algebra goes to an ordinary $gl(n)$ algebra when the deformation parameters
$q$ and $p$ goes to 1.

\def\ot{\otimes}
\section{q-symmetric states}
In this section we study the statistics of many particle state.
Let $N$ be the number of particles. Then the N-partcle state can be obtained
from
the tensor product of single particle state:
\bb
|i_1,\cdots,i_N>=|i_1>\ot |i_2>\ot \cdots \ot |i_N>,
\ee
where $i_1,\cdots, i_N$ take one value among $\{ 1,2,\cdots,n \}$ and the sigle
particle state is defined by $|i_k>=\a_{i_k}|0>$.

Consider the case that k appears $n_k$ times in the set $\{ i_1,\cdots,i_N\}$.
Then we have
\bb
n_1 + n_2 +\cdots + n_n =\sum_{k=1}^n n_k =N.
\ee
Using these facts we can define the q-symmetric states as follows:
\bb 
|i_1,\cdots, i_N>_q
=\sqrt{\frac{[n_1]_{p^2}!\cdots [n_n]_{p^2}!}{[N]_{p^2}!}}
\sum_{\sigma \in Perm}
\mbox{sgn}_q(\sigma)|i_{\sigma(1)}\cdots i_{\sigma(N)}>,
\ee
where
\begin{displaymath}
\mbox{sgn}_q(\sigma)=
q^{R(i_1\cdots i_N)}p^{R(\sigma(1)\cdots \sigma(N))},
\end{displaymath}
\begin{equation}
R(i_1,\cdots,i_N)=\sum_{k=1}^N\sum_{l=k+1}^N R(i_k,i_l)
\end{equation}
and $[x]_{p^2}=\fr{p^{2x}-1}{p^2-1}$.
Then the q-symmetric states obeys
\begin{equation}
|\cdots, i_k,i_{k+1},\cdots>_q=
\cases{
q^{-1} |\cdots,i_{k+1},i_k,\cdots>_q & if $i_k<i_{k+1}$\cr
 |\cdots,i_{k+1},i_k,\cdots>_q & if $i_k=i_{k+1}$\cr
q |\cdots,i_{k+1},i_k,\cdots>_q & if $i_k>i_{k+1}$\cr
}
\end{equation}
The above property can be rewritten by introducing the deformed transition
operator
$P_{k,k+1}$ obeying
\bb
P_{k,k+1}
|\cdots, i_k , i_{k+1},\cdots>_q =|\cdots, i_{k+1},i_k,\cdots>_q
\ee
This operator satisfies 
\bb
P_{k+1,k}P_{k,k+1}=Id,~~~\mbox{so}~~P_{k+1,k}=P^{-1}_{k,k+1}
\ee
Then the equation (33) can be written as
\bb
P_{k,k+1}
|\cdots, i_k , i_{k+1},\cdots>_q 
=q^{-\epsilon(i_k,i_{k+1})}
|\cdots, i_{k+1},i_k,\cdots>_q
\ee
where $\epsilon(i,j)$ is defined as
\bd
\epsilon(i,j)=
\cases{
1 & if $ i>j$\cr
0 & if $ i=j$ \cr
-1 & if $ i<j$ \cr }
\ed
It is worth noting that the relation (36) does not contain the deformation
parameter
 $p$. And the relation (36) goes to the symmetric relation for the ordinary
bosons 
 when the deformation parameter $q$ goes to $1$.
If we define the fundamental q-symmetric state $|q>$ as
\bd
|q>=|i_1,i_2,\cdots,i_N>_q
\ed
with $i_1 \leq i_2 \leq \cdots \leq i_N$, we have for any $k$
\bd
|P_{k,k+1}|q>|^2 =||q>|^2 =1.
\ed
In deriving the above relation we used following identity
\bd
\sum_{\sigma \in Perm } p^{R(\sigma(1),\cdots, \sigma(N))}=
\fr{[N]_{p^2}!}{[n_1]_{p^2}!\cdots [n_n]_{p^2}!}.
\ed

\section{Concluding Remark}
To conclude,  I used the two parameter deformed multimode oscillator system
given in ref [12] to construct its representation, coherent states and deformed
$gl_q(n)$ algebra.
Mutimode oscillator is important when we investigate the many body quantum
mechanics 
and statistical mechanics.
In order to construct the new statistical behavior for deformed particle
obeying the algebra (1), I investigated the defomed symmetric property of
two parameter deformed mutimode states.

\section*{Acknowledgement}
This                   paper                was
supported         by  
the   KOSEF (961-0201-004-2)   
and   the   present   studies    were   supported   by   Basic  
Science 
Research Program, Ministry of Education, 1995 (BSRI-95-2413).

\vfill\eject

\end{document}